\documentstyle[twoside]{article}

\catcode`\@=11
\long\def\@makefntext#1{
\protect\noindent \hbox to 3.2pt {\hskip-.9pt
$^{{\eightrm\@thefnmark}}$\hfil}#1\hfill}               
\def\thefootnote{\fnsymbol{footnote}}
\def\@makefnmark{\hbox to 0pt{$^{\@thefnmark}$\hss}}    
\def\ps@myheadings{\let\@mkboth\@gobbletwo
\def\@oddhead{\hbox{}
\rightmark\hfil\eightrm\thepage}
\def\@oddfoot{}\def\@evenhead{\eightrm\thepage\hfil
\leftmark\hbox{}}\def\@evenfoot{}
\def\sectionmark##1{}\def\subsectionmark##1{}}

\oddsidemargin=\evensidemargin
\renewcommand{\thefootnote}{\fnsymbol{footnote}}

\newcounter{sectionc}\newcounter{subsectionc}
\newcounter{subsubsectionc}
\renewcommand{\section}[1] 
{\vspace{12pt}\addtocounter{sectionc}{1}
\setcounter{subsectionc}{0}\setcounter{subsubsectionc}{0}
\noindent
        {\tenbf\thesectionc. #1}\par\vspace{5pt}}
\renewcommand{\subsection}[1] {\vspace{12pt}
\addtocounter{subsectionc}{1}
        \setcounter{subsubsectionc}{0}\noindent
        {\bf\thesectionc.\thesubsectionc. 
{\kern1pt \bfit #1}}\par\vspace{5pt}}
\renewcommand{\subsubsection}[1] 
{\vspace{12pt}\addtocounter{subsubsectionc}{1}
        \noindent{\tenrm\thesectionc.\thesubsectionc.
\thesubsubsectionc.
        {\kern1pt \tenit #1}}\par\vspace{5pt}}
\newcommand{\nonumsection}[1] {\vspace{12pt}\noindent{\tenbf #1}
        \par\vspace{5pt}}

\newcounter{appendixc}
\newcounter{subappendixc}[appendixc]
\newcounter{subsubappendixc}[subappendixc]
\renewcommand{\thesubappendixc}{\Alph{appendixc}.
\arabic{subappendixc}}
\renewcommand{\thesubsubappendixc}
        {\Alph{appendixc}.\arabic{subappendixc}.
\arabic{subsubappendixc}}

\renewcommand{\appendix}[1] {\vspace{12pt}
        \refstepcounter{appendixc}
        \setcounter{figure}{0}
        \setcounter{table}{0}
        \setcounter{lemma}{0}
        \setcounter{theorem}{0}
        \setcounter{corollary}{0}
        \setcounter{definition}{0}
        \setcounter{equation}{0}
        \renewcommand{\thefigure}{\Alph{appendixc}.
\arabic{figure}}
        \renewcommand{\thetable}{\Alph{appendixc}.
\arabic{table}}
        \renewcommand{\theappendixc}{\Alph{appendixc}}
        \renewcommand{\thelemma}{\Alph{appendixc}.
\arabic{lemma}}
        \renewcommand{\thetheorem}{\Alph{appendixc}.
\arabic{theorem}}
        \renewcommand{\thedefinition}{\Alph{appendixc}.
\arabic{definition}}
        \renewcommand{\thecorollary}{\Alph{appendixc}.
\arabic{corollary}}
        \renewcommand{\theequation}{\Alph{appendixc}.
\arabic{equation}}
        \noindent{\tenbf Appendix \theappendixc #1}
\par\vspace{5pt}}
\newcommand{\subappendix}[1] {\vspace{12pt}
        \refstepcounter{subappendixc}
        \noindent{\bf Appendix \thesubappendixc. 
{\kern1pt \bfit #1}}
        \par\vspace{5pt}}
\newcommand{\subsubappendix}[1] {\vspace{12pt}
        \refstepcounter{subsubappendixc}
        \noindent{\rm Appendix \thesubsubappendixc. 
{\kern1pt \tenit #1}}
        \par\vspace{5pt}}

\newcommand{\textlineskip}{\baselineskip=13pt}
\newcommand{\smalllineskip}{\baselineskip=10pt}

\def\eightcirc{
\begin{picture}(0,0)
\put(4.4,1.8){\circle{6.5}}
\end{picture}}
\def\eightcopyright{\eightcirc\kern2.7pt\hbox{\eightrm c}}

\def\abstracts#1#2#3{{
        \centering{\begin{minipage}{4.5in}
        \baselineskip=10pt\footnotesize
        \parindent=0pt #1\par
        \parindent=15pt #2\par
        \parindent=15pt #3
        \end{minipage}}\par}}

\renewenvironment{thebibliography}[1]                   
        {\frenchspacing
         \ninerm\baselineskip=11pt
         \begin{list}{\arabic{enumi}.}
        {\usecounter{enumi}\setlength{\parsep}{0pt}
         \setlength{\leftmargin 12.7pt}{\rightmargin 0pt} 
         \setlength{\itemsep}{0pt} \settowidth
        {\labelwidth}{#1.}\sloppy}}{\end{list}}

\newcounter{itemlistc}
\newcounter{romanlistc}
\newcounter{alphlistc}
\newcounter{arabiclistc}

\newcommand{\fcaption}[1]{
        \refstepcounter{figure}
        \setbox\@tempboxa = 
       \hbox{\footnotesize Fig.~\thefigure. #1}
        \ifdim \wd\@tempboxa > 5in
           {\begin{center}
        \parbox{5in}{\footnotesize\smalllineskip 
        Fig.~\thefigure. #1}
            \end{center}}
        \else
             {\begin{center}
             {\footnotesize Fig.~\thefigure. #1}
              \end{center}}
        \fi}

\newcommand{\tcaption}[1]{
        \refstepcounter{table}
        \setbox\@tempboxa = \hbox{\footnotesize 
          Table~\thetable. #1}
        \ifdim \wd\@tempboxa > 5in
           {\begin{center}
        \parbox{5in}{\footnotesize\smalllineskip 
         Table~\thetable. #1}
            \end{center}}
        \else
             {\begin{center}
             {\footnotesize Table~\thetable. #1}
              \end{center}}
        \fi}

\def\@citex[#1]#2{\if@filesw\immediate\write\@auxout
        {\string\citation{#2}}\fi
\def\@citea{}\@cite{\@for\@citeb:=#2\do
        {\@citea\def\@citea{,}\@ifundefined
        {b@\@citeb}{{\bf ?}\@warning
        {Citation `\@citeb' on page \thepage \space undefined}}
        {\csname b@\@citeb\endcsname}}}{#1}}

\newif\if@cghi
\def\cite{\@cghitrue\@ifnextchar [{\@tempswatrue
        \@citex}{\@tempswafalse\@citex[]}}
\def\citelow{\@cghifalse\@ifnextchar [{\@tempswatrue
        \@citex}{\@tempswafalse\@citex[]}}
\def\@cite#1#2{{$\null^{#1}$\if@tempswa\typeout
        {IJCGA warning: optional citation argument
        ignored: `#2'} \fi}}

\def\pmb#1{\setbox0=\hbox{#1}
        \kern-.025em\copy0\kern-\wd0
        \kern.05em\copy0\kern-\wd0
        \kern-.025em\raise.0433em\box0}

\def\fnt#1#2{\footnotetext{\kern-.3em
        {$^{\mbox{\scriptsize #1}}$}{#2}}}

\def\fpage#1{\begingroup
\voffset=.3in
\thispagestyle{empty}\begin{table}[b]
\centerline{\footnotesize #1}
        \end{table}\endgroup}

\def\runninghead#1#2{\pagestyle{myheadings}
\markboth{{\protect\footnotesize\it{\quad #1}}\hfill}
{\hfill{\protect\footnotesize\it{#2\quad}}}}
\font\tenrm=cmr10
\font\tenit=cmti10
\font\tenbf=cmbx10
\font\bfit=cmbxti10 at 10pt
\font\ninerm=cmr9
\font\nineit=cmti9
\font\ninebf=cmbx9
\font\eightrm=cmr8

\def\qed{\hbox{${\vcenter{\vbox{                        
   \hrule height 0.4pt\hbox{\vrule width 0.4pt height 6pt
   \kern5pt\vrule width 0.4pt}\hrule height 0.4pt}}}$}}

\renewcommand{\thefootnote}{\fnsymbol{footnote}}        

\def\bsc{{\sc a\kern-6.4pt\sc a\kern-6.4pt\sc a}}       
\def\bflatex{\bf L\kern-.30em\raise.3ex\hbox{\bsc}\kern-.14em
T\kern-.1667em\lower.7ex\hbox{E}\kern-.125em X}

\begin{document}

\def\be{\begin{equation}}  
\def\ee{\end{equation}} 
\def\ea{\end{array}\end{equation}}  
\def\bac{\begin{equation}\begin{array}{rll}}
\newcommand{\uq}{U_q (\widehat{sl(2)})}
\def\gam{\gamma}  
\def\la{\lambda}  
\def\ra{\rightarrow}  
\def\va{\varphi}  
\def\pa{\partial}  
\def\ps{\psi}
\def\Zp{{\Bbb Z}>0}  
\def\Zpn{{\Bbb Z}/\{ 0\}}  
\def\Zn{{\Bbb Z}<0}  
\def\Z{{\Bbb Z}}  
\def\C{{\Bbb C}}  
\def\N{{\Bbb N}} 
\def\ep{\epsilon}  
\def\epp{\epsilon^{\prime}}
\def\s{\sigma}
\newcommand{\eq}{\end{equation}}
\newcommand{\bq}{\begin{equation}}
\def\Zp{{\Bbb Z}>0}
\def\ov{\over}
\def\pl{\prod\limits}
\runninghead{A.H. Bougourzi, M. Couture and M. Kacir} 
{From quantum groups to the exact dynamical...}
\normalsize\textlineskip
\thispagestyle{empty}
\setcounter{page}{1}

\rightline{ITP-SB-96-28}  
\rightline{ June, 1996}
\vspace*{0.88truein}

\fpage{1}
\centerline{\bf FROM QUANTUM AFFINE GROUPS TO THE 
EXACT DYNAMICAL} 
\vspace*{0.035truein}
\centerline{\bf  CORRELATION FUNCTION OF THE 
HEISENBERG MODEL}
\vspace*{0.37truein}
\centerline{\footnotesize A.H. BOUGOURZI}
\vspace*{0.015truein}
\centerline{\footnotesize\it Institute of Theoretical Physics, 
SUNY at Stony Brook}
\baselineskip=10pt
\centerline{\footnotesize\it Stony Brook, NY 11794}
\vspace*{10pt}
\centerline{\footnotesize M. COUTURE}
\vspace*{0.015truein}
\centerline{\footnotesize\it Neutron \& Condensed Matter 
Sciences,  Chalk River Laboratories}
\baselineskip=10pt
\centerline{\footnotesize\it Chalk River, 
Ontario, Canada KOJ 1J0}
\vspace*{10pt}
\centerline{\normalsize and}
\vspace*{10pt}
\centerline{\footnotesize M. KACIR}
\vspace*{0.015truein}
\centerline{\footnotesize\it Service de Physique Theorique,
 CE-Saclay}
\baselineskip=10pt
\centerline{\footnotesize\it F-91191 Gif-sur-Yvette Cedex, 
France}
\vspace*{0.225truein}

\vspace*{0.21truein}
\abstracts{The exact form factors of the Heisenberg models
$XXX$ and $XXZ$ have been recently 
computed through the quantum affine symmetry of $XXZ$ model 
in the
thermodynamic limit. We use them to derive an exact formula
for the contribution of two spinons to the dynamical 
correlation function of $XXX$ model at zero temperature.}{}{}



\vspace*{1pt}\textlineskip  
\section{Introduction}    
\vspace*{-0.5pt}
\noindent
This letter is a short summary of our results on 
the two-spinon
dynamical correlation function (DCF) of the Heisenberg
model$^1$. 
Scattering cross sections of neutrons with spin chains are
directly proportional to dynamical correlation functions.
The latter, in turn, are expressed as series in terms of
form factors of  local spin operators. 
Unfortunately, unlike those in quantum field 
theories,$^2$ form factors of spin chains are
very complicated  because of the non-relativistic dispersion
relations of these lattice models. For more details on the
usefulness of form factors and DCF
see Refs.$^{3,4,5,6,7,8,9}$. In a recent development however, 
it has
been realized that in the thermodynamic limit the Heisenberg
model $XXZ$ becomes symmetric under the quantum affine
algebra $U_q(\widehat{sl(2)})$ in the anti-ferromagnetic 
regime.$^{10}$ 
This puts its resolution  on the same footing as that of
 conformal
field theory models. In fact, the bosonization technic of
conformal field theory extends
to this case and has been used to compute exact  static
correlation functions and form factors of spin local operators
of $XXZ$, and also those of $XXX$ after taking the isotropic 
limit $q\rightarrow -1$ of the former.$^{10}$ 
Unfortunately, the latter
physical quantities have very cumbersome multi-integral form
which limit their usefulness for the exact computation of 
DCF. However, our main point here
is that there is one exception to this latter
statement and that is the form factors needed for the
contribution of two spinons to the DCF
of the Heisenberg model  $XXX$ have more tractable  form.
Therefore we use them to compute this more interesting quantity 
of the exact two-spinon DCF of $XXX$ in the thermodynamic
limit and at zero temperature.

\textheight=7.8truein
\setcounter{footnote}{0}
\renewcommand{\thefootnote}{\alph{footnote}}

\section{Quantum affine symmetry of $XXZ$ model}
\noindent
In this section, we briefly review the $XXZ$ model in the 
thermodynamic limit and its quantum affine symmetry. 
This model is defined through its Hamiltonian
\be { H_{XXZ}}=
-{1\over 2} \sum_{n=-\infty}^{\infty} (\s_n^x \s_{n+1}^x 
+\s_n^y \s_{n+1}^y + \Delta\s_n^z \s_{n+1}^z) 
\label{hamiltonian}, 
\ee 
where $\Delta=(q+q^{-1})/2$ is the anisotropy parameter. 
Here $\sigma_n^{x,y,z}$ 
are the usual
Pauli matrices acting at the $n^{\rm th}$ position of the 
formal 
infinite tensor product 
\be W= \cdots V \otimes V \otimes V  \cdots 
\label{infprod},
\eq 
where $V$ is the two-dimensional representation of $U_q(sl(2))$ 
quantum group. 
We consider the model in 
the anti-ferromagnetic regime $\Delta <-1$, i.e., 
$-1<q<0$. The  action of $H_{XXZ}$ 
on $W$ is
not well defined due to the appearance of 
divergences. However, since this model is symmetric 
under the quantum affine group $\uq$,  the eigenspace on which
this action becomes well defined 
is identified with the following level 0 $\uq$ 
module:
\be
{\cal F}=\sum_{i,j}V(\Lambda_i)\otimes V(\Lambda_j)^*, 
\ee
where $\Lambda_i$ and $V(\Lambda_i); i=0,1$ are  level 
1 $\uq$-highest
weights and $\uq$-highest weight modules, respectively. 
Roughly speaking, $V(\Lambda_i)$ is identified 
with the
subspace of the formal semi-infinite space
\be X= \cdots V \otimes V \otimes V,\eq 
consisting  of all linear combinations of  
spin configurations with  fixed boundary conditions such that 
the eigenvalues of $\sigma^z_n$ are $(-1)^{i+n}$ in the limit 
$n\ra -\infty$. The eigenspace ${\cal F}$ consists of spinon 
particles$^{11}$
$\{|\xi_1,\cdots 
\xi_n>_{\ep_1,\cdots \ep_n;i}, n\geq 0\}$. Here $i$ fixes
the boundary condition, $\xi_j$ are the spectral parameters,
and $\ep_j=\pm 1$ are 
the spins of the
spinons. The completeness relation reads:$^{10}$ 
\be {\bf I}=\sum_{i=0,1}\sum_{n \geq 0} \sum_{\ep_1,\cdots,
\ep_n=\pm 1}
{1 \ov {n !}} \oint  {d\xi_1\over 2\pi i \xi_1} \cdots 
 {d\xi_n\over 2\pi i \xi_n}
 |\xi_n,\cdots,\xi_1>_{{\ep_n,\cdots,\ep_1};i}\;
{_{i;{\ep_1,\cdots,\ep_n}}{<\xi_1,\cdots,\xi_n|}}. 
\ee

The actions of $H_{XXZ}$ and the translation operator 
$T$, which shifts the spin  chain by one site, on ${\cal F}$ 
are given by 
\bac T|\xi_1,\cdots,\xi_n>_i &=&\pl_{i=1}^n\tau(\xi_i)^{-1}
|\xi_1,\cdots,\xi_n>_{1-i},\quad 
T|0>_i =|0>_{1-i}, \\
H_{XXZ}|\xi_1,\cdots,\xi_n>_i 
&=&\sum_{i=1}^n e(\xi_i)|\xi_1,\cdots,\xi_n>_i, 
\label{states}\ea  
where
\bac \tau(\xi)&=& \xi^{-1} 
{\theta_{q^4} (q \xi^2) \ov \theta_{q^4} (q 
\xi^{-2})}=e^{-i p(\alpha)},\quad p(\alpha)=am({2K\over \pi}
\alpha)-{\pi/2}, \\
e(\xi) &=&{ 1-q^2 \ov 2 q} \xi {d \ov d \xi} 
\log \tau(\xi)= {2K\over \pi} \sinh({\pi K^\prime\over K})
dn({2K\over \pi}\alpha). \label{enmom}
\ea
Here, $e(\xi)$ and  $p(\alpha)$ are the energy and the 
momentum of
the spinon respectively,  $am(x)$ and $dn(x)$ are the usual 
elliptic amplitude and
delta functions, with the complete elliptic integrals
 $K$ and $K^\prime$, and
\bac
q&=&-\exp(-\pi K^\prime/K),\\
\xi&=&ie^{i\alpha},\\ 
\theta_x(y)&=&(x;x)_{\infty} (y;x)_{\infty} 
(x y^{-1};x)_{\infty},\\
(y;x)_{\infty}&=&\prod_{n=0}^{\infty} (1-y x^n).
\ea
This means, $\sigma^{x,y,z}(t, n)$ at time $t$ and position
$n$ are related to $\sigma^{x,y,z}(0,0)$ at time 0 and position
0 through:
\be \sigma^{x,y,z} (t,n)=\exp(i tH_{XXZ} ) T^{-n} 
\sigma^{x,y,z} 
(0,0) T^{n} 
\exp(-i tH_{XXZ} ). \ee

\section{Two-spinon dynamical correlation function of $XXX$
model}
\noindent

Here we first define one of the components of the DCF in the
case of the 
$XXZ$ model , where the spinon picture is well understood,
and then take its isotropic limit to get {\it all} 
non-vanishing components of the DCF of $XXX$  model using its
isotropy. We consider the component
\be
S^{i,+-}(w, k)=
\int_{-\infty}^{\infty} dt \sum_{n\in Z}
e^{i(wt+kn)} {_i}< 0|\sigma^+(t, n)\sigma^-(0,0)|0>_i,
\ee
where $w$ and $k$ are the neutron energy and momentum transfer
respectively, and $i$ corresponds to the boundary condition. 
Later, we 
find that the DCF is in fact independent of $i$.
Using the completeness relation, 
the two-spinon contribution is given by
\bac
&&S_2^{i,+-}(w,k)= {\pi} 
\sum_{n\in Z} \sum_{\ep_1,\ep_2} 
\oint {d\xi_1\over 2\pi i \xi_1}
 {d\xi_2\over 2\pi i \xi_2}
\exp\left(in(k+p(\xi_1)+p(\xi_2))\right)
\\
&&\times \delta (w-e(\xi_1)
-e(\xi_2)){_{i+n}<0|}\sigma^+ 
(0,0)|\xi_2,\xi_1>_{\ep_2,\ep_1;i+n}\:
{_{i;\ep_1,\ep_2}<\xi_1,\xi_2|}
\sigma^-(0,0)|0>_i.
\ea
It can be put in the following tractable form:
\bac
&& S_2^{i,+-}(w,k)=\pi
\sum_{\ep_1,\ep_2} 
 \oint {d \xi_1\over 2\pi i \xi_1}
{ d \xi_2\over 2 \pi i \xi_2} 
  \sum_{n\in Z} \exp\left( 2 i n(k+p(\xi_1)+p(\xi_2))\right)
\\
&&\times \delta(w-e(\xi_1)-e(\xi_2)) \Bigl(
{_{i}<0|}\sigma^+(0,0)|\xi_2,\xi_1>_{\ep_2,\ep_1;i}
\:
{_{i;\ep_1,\ep_2}<\xi_1,\xi_2|}
\sigma^-{(0,0)}|0>_i \Bigr.\\
&& +
\left. \exp\left( i(k+p(\xi_1)+p(\xi_2)\right)
{_{1-i}<0|}\sigma^+(0,0)|\xi_2,\xi_1>_{\ep_2,
\ep_1;1-i}\right.\\
&& \:
{_{i;\ep_1,\ep_2}<\xi_1,\xi_2|}
\sigma^-{(0,0)}|0>_i \Bigr).
\label{corr2}\ea
Now the various form factors invloved in this expression have
been computed in Ref. 10. For our purposes, we  give
only their isotropic limits which are obtained by first making
the redefinitions
\bac
\xi&=&ie^{{\epsilon\beta\over i\pi}},\\
q&=&-e^{-\epsilon},\quad\epsilon\ra 0^+,
\ea 
with $\beta$ being the appropriate 
spectral parameter  for $XXX$ model, and then taking the limit
$q\ra -1$. Therefore one finds:$^{10}$ 
\bac
&&|{_{i}{<0|}}\sigma^+(0,0)|\xi_2,\xi_1>_{--;i}|^2
{d\xi_{1}\over 2\pi i  \xi_{1}} 
{d\xi_{2}\over  2\pi i \xi_{2}}\ra
{ \Gamma(3/4)^2|A_{-}(\beta_1-\beta_2)|^2 d\beta_1 d\beta_2\over
16\Gamma(1/4)^2 |A_+(i\pi/2)|^2 |A_{-}(i\pi/2)|^2
\cosh(\beta_1)\cosh(\beta_2)},\\
&&p(\xi)\ra p(\beta),\quad {\rm
s.t.}\quad \cot(p(\beta))=\sinh(\beta),\quad -\pi\leq p(\beta)
\leq 0,\\
&&e(\xi)\ra e(\beta)={\pi\over \cosh(\beta)}=-\pi \sin(p(\beta)),
\ea
where
\be
|A_{\pm}(\beta)|^2=
\exp\left( -\int_{0}^{\infty}dx {(\cosh(2 x(1-{\delta\over\pi}))
\cos({2 x \gamma
\over \pi})-1)
\exp(\mp x)\over x \sinh(2 x)\cosh(x)}\right).
\ee
Here $\Gamma(x)$ is the usual gamma function and 
$\beta=\gamma+i\delta$, with $\gamma$ and $\delta$ being real.
Restricting
to the first Brillouin zone, integrating  the continuous and
discrete delta functions and keeping track of the Jacobian 
factors, we find that the latter expression is independent
of the boundary conditions $i$ (which is henceforth 
omitted). Moreover, it substantially
simplifies to:
\be
S^{+-}_2(w, k)={\pi^2 {\Gamma(3/4)^2}
\Theta(2\pi \sin(k/2)-w)\Theta(
w-\pi|\sin(k)|)\over 4 \Gamma(1/4)^2
|A_{-}(i\pi/2)|^2|A_{+}(i\pi/2)|^2}
{|A_{-}(\bar\beta_1-\bar\beta_2)|^2\over
\sqrt{(2\pi \sin(k/2))^2-w^2}},
\label{exact}
\ee
where $\Theta$ is the Heaviside step function, and 
for fixed $w$ and $k$,  
$\bar\beta_1$ and $\bar\beta_2$ are the solutions to:
\bac
w&=&e(\bar\beta_1)+e(\bar\beta_2),\\
k&=&-p(\bar\beta_1)-p(\bar\beta_2).
\ea
Note that the pair $(\bar\beta_1, \bar\beta_2)$ 
is identified with the pair
$(\bar\beta_2, \bar\beta_1)$.
From the isotropy of the Heisenberg model and the inclusion 
of both sectors $i=0$ and $i=1$,
we  obtain all   non-vanishing components of its DCF from
$S^{+-}_2(w,k)$ as:
\be
S^{xx}_2(w,k)=S^{yy}_2(w,k)=S^{zz}_2(k,w)=4S^{+-}_2(w,k),
\ee 
with 
\be
\sigma^{\pm}={\sigma^x\pm i\sigma^y\over 2}.
\ee 
Despite its square root singularity, $S_2^{+-}(w,k)$
actually vanishes in the vicinity of the upper boundary 
$w_u=2\pi\sin(k/2)$ in the dispersion relation  of two spinons.
Moreover, it diverges in the vicinity of the lower boundary
which is given by the des Cloizeaux-Pearson dispersion 
relation $w_l=\pi|\sin k|$. This behaviour reproduces very
well the one obtained previously through the ansatz made
in Ref. 5, although now the upper cutoff appears naturally.
It would be interesting to
investigate to which extent the two-spinon DCF 
verifies the sum rules which are valid for the full DCF.$^{12}$
The extension of
this work to the Heisenberg  model with higher spin is certainly
desirable. In this case, the form factors can in 
principle be computed through the bosonization of the vertex
operators which is now available in Ref. 13.

\nonumsection{Acknowledgements}
\noindent
The work of A.H.B. is supported by the NSF Grant \#
PHY9309888.
We wish to thank Elgradechi, Karbach, Korepin, Lorenzano, 
McCoy, M\"uller, Perk, Shrock, Sebbar,
Takhtajan and Weston for useful discussions.

\nonumsection{References}

\end{document}